\documentclass[a4paper,fleqn]{cas-sc}

\usepackage[numbers,sort&compress]{natbib}
\usepackage{amsmath,amssymb,amsfonts}
\usepackage{mathrsfs}
\usepackage{placeins}
\usepackage{float}
\usepackage{multirow}
\usepackage{algorithm}
\usepackage{algorithmicx}
\usepackage{algpseudocode}
\usepackage{listings}
\usepackage{booktabs}

\newtheorem{proposition}{Proposition}

\makeatletter
\newenvironment{inlinefigure}{%
  \par\medskip\noindent\begin{minipage}{\linewidth}\centering\def\@captype{figure}%
}{%
  \end{minipage}\par\medskip
}
\newenvironment{inlinetable}{%
  \par\medskip\noindent\begin{minipage}{\linewidth}\centering\def\@captype{table}%
}{%
  \end{minipage}\par\medskip
}
\makeatother

\begin{document}
\let\WriteBookmarks\relax
\def\floatpagepagefraction{1}
\def\textpagefraction{.001}

\shorttitle{Recentered-Domain Yau--Yau Filter for Target Tracking}
\shortauthors{L. Ma et~al.}

\title[mode=title]{A Recentered-Domain Yau--Yau Filter for Target Tracking}

\author[1,2,3]{Lei Ma}
\credit{Conceptualization, Methodology, Software, Validation, Formal analysis, Investigation, Data curation, Visualization, Writing - original draft}

\author[1,2,3]{Yuzhong Hu}
\credit{Methodology, Formal analysis, Investigation, Writing - review and editing}

\author[3]{Xiaoming Zhang}
\cormark[1]
\ead{zhangxiaoming@bimsa.cn}
\credit{Supervision, Methodology, Project administration, Writing - review and editing}

\affiliation[1]{organization={University of Chinese Academy of Sciences},
                city={Beijing},
                postcode={100049},
                country={China}}

\affiliation[2]{organization={Academy of Mathematics and Systems Science, Chinese Academy of Sciences},
                city={Beijing},
                postcode={100190},
                country={China}}

\affiliation[3]{organization={Beijing Institute of Mathematical Sciences and Applications},
                city={Beijing},
                postcode={101408},
                country={China}}

\cortext[cor1]{Corresponding author.}

\begin{abstract}
The Yau--Yau filter reformulates nonlinear state estimation as probability-density propagation governed by the Forward Kolmogorov equation (FKE). Applying it to target tracking, however, requires efficient FKE approximation on a finite computational domain. This paper proposes a Recentered-Domain Yau--Yau Filter (RD-YYF), which solves the FKE within a fixed-size local window centered at the latest state estimate. This design concentrates numerical resolution near the dominant posterior density. Offline, physics-informed neural networks (PINNs) generate FKE solution snapshots, while principal component analysis constructs a low-dimensional representation of density evolution. A lightweight residual surrogate maps the initial-condition coefficients and domain center to the terminal-solution coefficients. Online, the pretrained surrogate predicts density evolution within the recentered window, followed by observation update and state estimation. Experiments on two geometrically constrained target-tracking examples show that RD-YYF achieves lower tracking errors than the extended Kalman filter (EKF) and particle filter (PF), while retaining efficient per-timestep inference. Ablation results indicate that domain recentering improves density approximation in high-probability regions and accelerates offline PINN convergence. These results demonstrate the potential of RD-YYF for efficient nonlinear target tracking.
\end{abstract}

\maketitle

\section{Introduction}\label{sec:introduction}

Target tracking is a fundamental problem in many engineering applications, including military surveillance \cite{ref1}, robotics \cite{ref2}, and autonomous navigation \cite{ref3}. Its main objective is to sequentially estimate the motion state of a target from noisy sensor observations. Filtering algorithms play a central role in this task because they provide recursive state estimates under uncertainty and allow tracking systems to operate in dynamic and noisy environments.

Classical filtering methods, including the Kalman filter \cite{ref4}, its nonlinear variants such as the extended Kalman filter (EKF) \cite{ref5,ref6}, and sampling-based methods such as the particle filter (PF) \cite{ref7}, have provided a solid foundation for target tracking. However, many practical tracking systems involve nonlinear dynamics, nonlinear observation functions, and maneuvering target motion, for which linear-Gaussian assumptions or finite-particle approximations may be insufficient. Recent studies have therefore introduced adaptive, robust, and application-specific filtering strategies for complex tracking conditions \cite{ref8,ref9,ref10,ref11}. Although these methods improve tracking performance in specific settings, they often depend on local approximations, distributional assumptions, or carefully tuned parameters, which may limit their robustness in strongly nonlinear systems.

In recent years, learning-based filtering methods have received increasing attention. Neural networks have been incorporated into filtering frameworks to learn nonlinear state-transition or observation relationships from data \cite{ref12}. For instance, Song et al. \cite{ref13} combined recurrent neural networks with Kalman filtering to improve nonlinear radar target tracking, and Shen et al. \cite{ref14} introduced a self-attention-based Transformer model for nonlinear maneuvering target tracking. In addition, application-specific filtering methods have been developed for complex environments such as underwater tracking \cite{ref15}. These approaches can achieve strong empirical performance when sufficient training data are available. Nevertheless, purely data-driven or highly task-specific methods may suffer from limited interpretability, sensitivity to distributional changes, and reduced transferability across different tracking scenarios.

An alternative route is provided by the Yau--Yau filtering framework \cite{ref16,ref17}, which reformulates nonlinear filtering through the evolution of the conditional probability density. Unlike methods based on linearization or Gaussian approximation, the Yau--Yau filter is derived from the Duncan--Mortensen--Zakai equation \cite{ref18,ref19,ref20} and can handle general nonlinear systems with arbitrary initial distributions. Its key advantage is that density propagation can be transformed into a Forward Kolmogorov equation (FKE) whose propagation step is independent of the observation path. This observation-independent structure is useful for target tracking because it separates density prediction from measurement update and supports an offline-online computational framework.

Despite its theoretical advantages, the practical implementation of the Yau--Yau filter remains challenging because it requires efficient numerical solution of the FKE. Existing numerical methods for Yau-type filtering problems, such as direct methods and Gaussian-approximation-based approaches \cite{ref21,ref22}, have improved the computational feasibility of solving the underlying Kolmogorov equations.
These methods reduce part of the computational burden, but they do not fully address the spatial-localization issue that arises when the dominant posterior density moves through the state space.
However, when the probability density must be approximated over a large finite computational domain, the resulting PDE computation can still be expensive. Moreover, in target tracking problems, the posterior density typically moves with the target state. A fixed computational domain may allocate substantial numerical resolution to low-probability regions, while reducing the effective approximation accuracy around the dominant posterior mass.

Physics-informed neural networks (PINNs) provide a flexible tool for solving PDEs by embedding governing equations into the training loss \cite{ref23}. This makes them suitable for generating FKE solutions without relying on a fixed mesh. However, directly deploying a PINN solver at every filtering step is still computationally prohibitive for real-time target tracking. Therefore, an efficient implementation of the Yau--Yau filter requires not only a PDE solver, but also a reduced-order representation and a localized computational strategy that concentrates approximation capacity on the high-probability region of the posterior density.

Motivated by these considerations, this paper proposes a Recentered-Domain Yau--Yau Filter (RD-YYF) for target tracking. To the best of our knowledge, this is the first study to adapt the Yau--Yau filtering framework to target-tracking problems. The proposed method addresses the main computational difficulty of applying Yau--Yau filtering in this setting, namely the efficient approximation of the FKE on a finite domain whose high-probability region moves with the target state.

The main contributions of this paper are summarized as follows. First, we introduce a recentered-domain Yau--Yau filtering formulation for geometrically constrained target tracking, which improves the effective numerical resolution around the dominant posterior density. Second, we construct an offline-online reduced-order FKE solver by combining PINN-generated solution snapshots, PCA-based density representation, and a lightweight residual surrogate network. Third, we evaluate RD-YYF on two nonlinear target-tracking examples and compare it with EKF and PF baselines, with an ablation study demonstrating the role of the recentered-domain strategy.

The remainder of this paper is organized as follows. Section~\ref{sec:preliminaries} reviews the mathematical background of continuous-time nonlinear filtering, the Yau--Yau filtering framework, PINNs, and PCA. Section~\ref{sec:framework} presents the proposed RD-YYF method, including the problem formulation, recentered-domain strategy, and offline-online execution procedure. Section~\ref{sec:simulations} reports the numerical experiments and ablation study. Section~\ref{sec:conclusion} concludes the paper and discusses future research directions.

\section{Preliminaries}\label{sec:preliminaries}

\subsection{Yau--Yau filter}\label{subsec:yauyau}

The Yau--Yau filter considers the type of continuous filtering problem described by the following signal-observation stochastic differential equations (SDE):
\begin{equation}
\left\{
\begin{aligned}
d x_t &= f(x_t,t)\,dt + G(x_t,t)\,d w_t,\\
d y_t &= h(x_t,t)\,dt + d v_t,\qquad y_0=0.
\end{aligned}
\right.
\label{eq:sde}
\end{equation}

Here, $t\in[0,T]$ denotes time, $x_t\in\mathbb{R}^n$ is the state variable, $y_t\in\mathbb{R}^m$ is the observation, $f$ and $h$ are $n$- and $m$-dimensional vector functions, respectively, and $G$ is an $n\times r$ matrix. The Brownian motions $w_t$ and $v_t$ satisfy
\begin{equation}
E[d w_t d w_t^\top]=Q_t\,dt,\qquad Q_t\in\mathbb{R}^{r\times r},
\qquad
E[d v_t d v_t^\top]=S_t\,dt,\qquad S_t\in\mathbb{R}^{m\times m},
\label{eq:brownian_cov}
\end{equation}
and are independent of each other and of the initial state $x_0$.

Let $\mathcal{F}_t=\sigma\{y_s,0\leq s\leq t\}$ denote the smallest $\sigma$-field generated by historical observations. Denote by $p(x,t)$ the conditional probability density of the state $x_t$ given historical observations. This conditional density can be obtained through the unnormalized density $\sigma(x,t)$ by solving the Duncan--Mortensen--Zakai (DMZ) equation \cite{ref18,ref19,ref20},
\begin{equation}
\left\{
\begin{aligned}
d\sigma(x,t) &= L\sigma(x,t)\,dt+\sigma(x,t) h^\top(x,t)S_t^{-1}\,d y_t,\\
\sigma(x,0)&=\sigma_0(x),
\end{aligned}
\right.
\label{eq:dmz}
\end{equation}
where
\begin{equation}
L(\cdot)=\frac{1}{2}\sum_{i,j=1}^n
\frac{\partial^2}{\partial x_i\partial x_j}
\left((GQG^\top)_{ij}\cdot\right)
-\sum_{i=1}^n\frac{\partial}{\partial x_i}\left(f_i\cdot\right).
\label{eq:operator}
\end{equation}

First, for each given observation path, the invertible exponential transformation \cite{ref24}
\begin{equation}
\sigma(x,t)=\exp\!\left(h^\top(x,t)S_t^{-1}y_t\right)\rho(x,t)
\label{eq:first_transform}
\end{equation}
transforms Eq.~\eqref{eq:dmz} into a deterministic partial differential equation with stochastic coefficients, which is referred to as the pathwise-robust DMZ equation:
\begin{equation}
\left\{
\begin{aligned}
\frac{\partial \rho}{\partial t}(x,t)
&+
\left[
\frac{\partial}{\partial t}
\left(h^\top(x,t)S_t^{-1}\right)y_t
\right]\rho(x,t) \\
&=
\exp\!\left(-h^\top(x,t)S_t^{-1}y_t\right)
\left(L-\frac{1}{2}h^\top S_t^{-1}h\right)
\left[
\exp\!\left(h^\top(x,t)S_t^{-1}y_t\right)\rho(x,t)
\right],\\
\rho(x,0)&=\sigma_0(x).
\end{aligned}
\right.
\label{eq:robust_dmz}
\end{equation}

Then, the observations arrive at discrete times $\tau_i=i\Delta t$, $i=0,1,2,\ldots,N_T$, with $\Delta t=T/N_T$. Let $\rho_i$ denote the solution of the robust DMZ equation with $y_t=y_{\tau_{i-1}}$ on the interval $t\in[\tau_{i-1},\tau_i]$, $i=1,2,\ldots,N_T$, namely
\begin{equation}
\left\{
\begin{aligned}
\frac{\partial \rho_i}{\partial t}(x,t)
&+
\left[
\frac{\partial}{\partial t}
\left(h^\top(x,t)S_t^{-1}\right)y_{\tau_{i-1}}
\right]\rho_i(x,t) \\
&=
\exp\!\left(-h^\top(x,t)S_t^{-1}y_{\tau_{i-1}}\right) 
\left(L-\frac{1}{2}h^\top S_t^{-1}h\right) 
\left[
\exp\!\left(h^\top(x,t)S_t^{-1}y_{\tau_{i-1}}\right)
\rho_i(x,t)
\right],\\
\rho_1(x,0)&=\sigma_0(x),\\
\rho_i(x,\tau_{i-1})
&=\rho_{i-1}(x,\tau_{i-1}),
\qquad i=2,3,\ldots,N_T.
\end{aligned}
\right.
\label{eq:piecewise_dmz}
\end{equation}

It is proved in \cite{ref17} that, both in the pointwise and $L^2$ senses,
\begin{equation}
\lim_{\Delta t\to0}
\sum_{i=1}^{N_T}
\mathbf{1}_{[\tau_{i-1},\tau_i)}(t)\rho_i(x,t)
=\rho(x,t).
\label{eq:rho_convergence}
\end{equation}

Therefore, $\rho_i(x,t)$ provides a good approximation to $\rho(x,t)$ on each observation interval.

With another exponential transformation, Eq.~\eqref{eq:piecewise_dmz} can be further transformed into an FKE with deterministic coefficients and stochastic initial values, as stated in the following proposition.

\begin{proposition}[\cite{ref17}]
For each $\tau_{i-1}\leq t\leq \tau_i$, $\rho_i(x,t)$ satisfies the robust DMZ equation if and only if
\begin{equation}
u_i(x,t)=\exp\!\left(h^\top(x,t)S_t^{-1}y_{\tau_{i-1}}\right)\rho_i(x,t)
\label{eq:second_transform}
\end{equation}
satisfies
\begin{equation}
\frac{\partial u_i}{\partial t}(x,t)
=\left(L-\frac{1}{2}h^\top S_t^{-1}h\right)u_i(x,t).
\label{eq:fke_general}
\end{equation}
\end{proposition}

Since the FKE in Eq.~\eqref{eq:fke_general} is independent of observations, it can be solved offline for a given initial condition. The entire procedure of the Yau--Yau filter consists of two steps.

\textit{Offline stage} $u_i(x,\tau_{i-1})\to u_i(x,\tau_i)$: solve the FKE on the time interval $[\tau_{i-1},\tau_i]$ with initial value $u_i(x,\tau_{i-1})$. This stage can be computed offline for any given initial values because the FKE is independent of the observations.

\textit{Online stage} $u_i(x,\tau_i)\to u_{i+1}(x,\tau_i)$: when the new observation $y_{\tau_i}$ arrives at time $\tau_i$, the initial value of $u_{i+1}(x,t)$ for the interval $[\tau_i,\tau_{i+1}]$ is updated by
\begin{equation}
\left\{
\begin{aligned}
u_1(x,0) &= \sigma_0(x),
&& t\in[0,\tau_1],\\
u_i(x,\tau_{i-1})
&= \exp\!\left[
h^\top(x,\tau_{i-1})S^{-1}_{\tau_{i-1}}
\left(y_{\tau_{i-1}}-y_{\tau_{i-2}}\right)
\right]   u_{i-1}(x,\tau_{i-1}),
&& t\in[\tau_{i-1},\tau_i],\ i\geq2.
\end{aligned}
\right.
\label{eq:yy_update}
\end{equation}

With these two computational stages over all time intervals, the functions $u_{i+1}(x,\tau_i)$ provide good approximations to the solution $\sigma(x,\tau_i)$ of the original DMZ equation at time $t=\tau_i$.

\subsection{Physics-informed neural networks}\label{subsec:pinn}

PINNs provide a mesh-free neural representation for approximating PDE solutions by embedding the governing equation and auxiliary conditions into the training loss. Compared with conventional finite difference or finite element schemes, PINNs avoid explicit mesh generation and can represent the solution as a continuous function of the input variables.

Unlike purely data-driven deep learning models that require large labeled datasets, PINNs embed the governing equations directly into the loss function as optimization constraints \cite{ref23,ref25}. By penalizing the PDE residual during training, the network is encouraged to satisfy the underlying system dynamics without relying on a prescribed spatial discretization.

The general structure of PINN is illustrated in Fig.~\ref{fig:pinn}. The network takes spatiotemporal coordinates as input and outputs an approximate solution to the given problem. The loss terms, which are computed from the network output and its partial derivatives, are efficiently evaluated using automatic differentiation. The total loss is formulated as a weighted combination of the residuals of the governing equation, the initial condition, and the boundary condition. By minimizing the total loss, the network parameters are optimized to ensure that the solution not only aligns with the training constraints but also adheres to the underlying physical laws, providing both accuracy and consistency.

\begin{inlinefigure}
\centering
\includegraphics[width=0.8\textwidth]{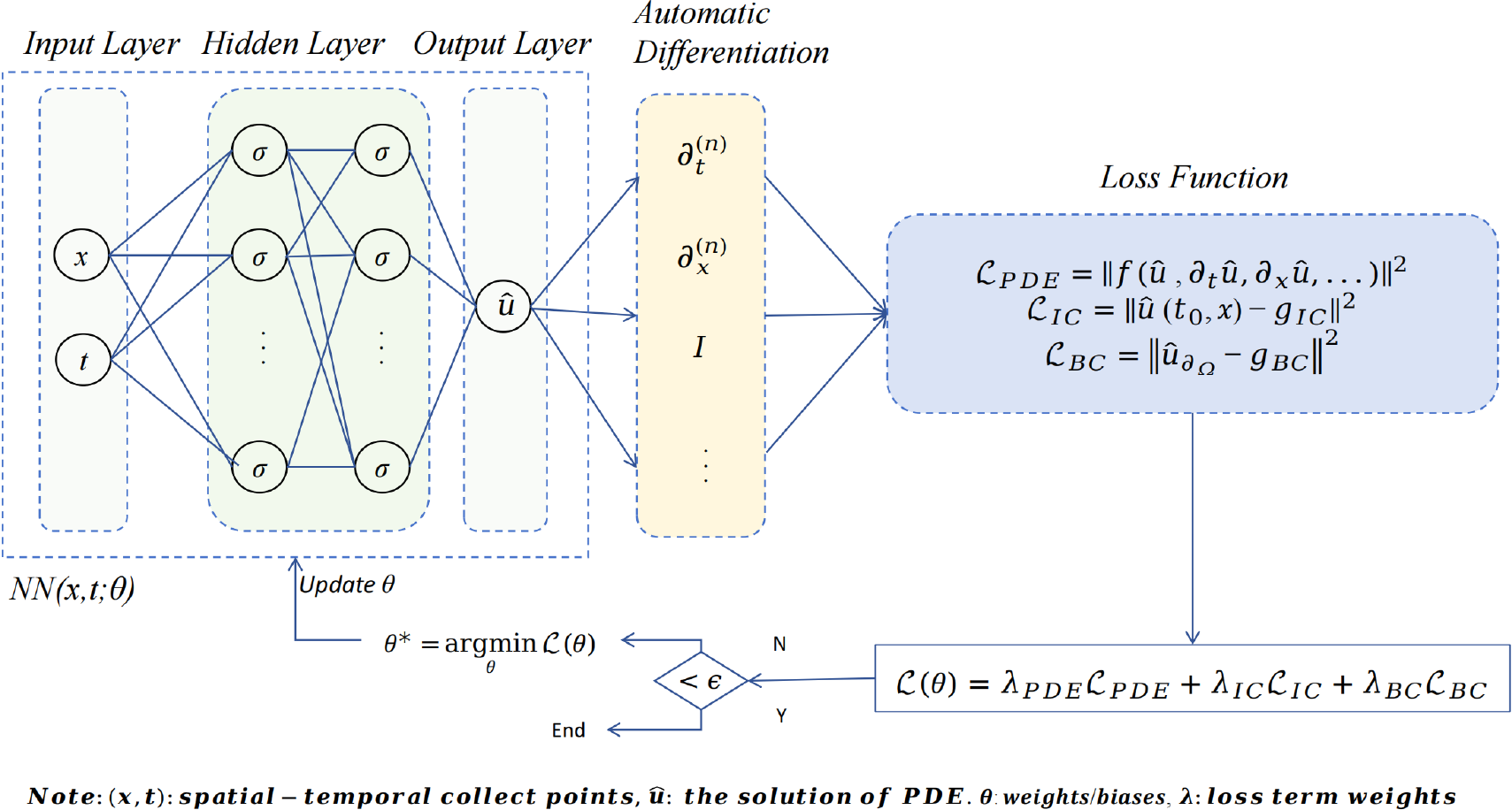}
\caption{The general structure of PINN.}
\label{fig:pinn}
\end{inlinefigure}

\subsection{Principal component analysis}\label{subsec:pca}

For high-dimensional solution snapshots or data matrices, dimensionality reduction is often needed to obtain compact representations while retaining the dominant modes of variation.

PCA is a widely used method for dimensionality reduction and feature extraction, particularly effective in high-dimensional settings. It transforms the original dataset into a new set of orthogonal components, called principal components (PCs), which are linear combinations of the original features \cite{ref26}. PCA begins by calculating the covariance matrix of the data and then derives its eigenvalues and eigenvectors. The eigenvalues quantify the variance captured by the corresponding eigenvectors, which define the PCs. These components serve as new axes along which the data are projected, with the first PC capturing the largest variance, followed sequentially by components with decreasing variance. By selecting a subset of these components, PCA effectively reduces dimensionality while preserving the most significant variability within the data \cite{ref27}.

\section{Proposed Framework}\label{sec:framework}

Based on the Yau--Yau filter, PINN, and PCA introduced in Section~\ref{sec:preliminaries}, this section constructs the proposed Recentered-Domain Yau--Yau Filter (RD-YYF). We first formulate the kinematically constrained target-tracking problem and derive its governing FKE in Section~\ref{subsec:problem}. We then introduce the recentered-domain strategy in Section~\ref{subsec:recentered} to localize the FKE computation around the dominant posterior region. Finally, we present the offline-online decoupled workflow in Section~\ref{subsec:offline_online} for efficient online state inference.

\subsection{Problem formulation and FKE derivation}\label{subsec:problem}

In many practical tracking scenarios, such as vehicles moving along a specific road or rail-based transportation systems, the target is constrained to follow a known geometric trajectory. Rather than artificially discretizing the system dynamics, we formulate the target-tracking problem in its natural continuous-time differential form. This formulation is well suited to the Yau--Yau filtering framework introduced in Section~\ref{sec:preliminaries}.

\subsubsection{Target-tracking problem definition}

Assume that the known trajectory curve can be parametrically represented as
\begin{equation}
x=\phi_x(s),\qquad y=\phi_y(s),
\label{eq:curve}
\end{equation}
where $s(t)$ is the arc-length parameter of the target along the trajectory, and $\phi_x(s)$ and $\phi_y(s)$ are known smooth functions describing the curve shape. With this geometric constraint, the target state can be uniquely determined by its arc length $s(t)$ and its tangential velocity $v(t)$ along the curve. We define the system state as $(s(t),v(t))$.

\textbf{\textit{State equation.}}
The target motion on this fixed trajectory satisfies the following continuous-time stochastic differential equations:
\begin{equation}
\left\{
\begin{aligned}
d s(t)&=v(t)\,dt+d\nu_s(t),\\
d v(t)&=f_v(s,t)\,dt+d\nu_v(t),
\end{aligned}
\right.
\label{eq:state_sv}
\end{equation}
where $f_v(\cdot)$ defines the tangential dynamics. The noise processes $\nu_s(t)$ and $\nu_v(t)$ are independent Wiener processes with covariance
\begin{equation}
\nu_t=\begin{bmatrix}\nu_s(t)\\ \nu_v(t)\end{bmatrix},\qquad
E[d\nu_t d\nu_t^\top]=Q_t\,dt.
\label{eq:state_noise}
\end{equation}

\textbf{\textit{Observation equation.}}
Although the target state evolves in the arc-length space $s(t)$, the observation data are typically acquired in the two-dimensional physical space. Suppose that the observation sensor is fixed at coordinates $(x_0,y_0)$ and provides range and bearing information. Even though measurements in practical target-tracking systems are usually collected at discrete time instants, it is common in continuous-time filtering to model the observation process in differential form. Specifically, we introduce the range observation process $R(t)$ and the bearing observation process $\Theta(t)$. These processes describe the evolution of observation information over time rather than instantaneous physical measurements:
\begin{equation}
\left\{
\begin{aligned}
dR(t)&=h_R(s)\,dt+d w_R(t),\\
d\Theta(t)&=h_\Theta(s)\,dt+d w_\Theta(t),
\end{aligned}
\right.
\label{eq:obs}
\end{equation}
where
\begin{equation} 
h_R(s)=\sqrt{(\phi_x(s)-x_0)^2+(\phi_y(s)-y_0)^2},\qquad
h_\Theta(s)=\arctan\frac{\phi_x(s)-x_0}{\phi_y(s)-y_0}.
\label{eq:obs_funcs}
\end{equation}

The observation noise satisfies
\begin{equation}
w_t=\begin{bmatrix}w_R(t)\\ w_\Theta(t)\end{bmatrix},\qquad
E[d w_t d w_t^\top]=S_t\,dt.
\label{eq:obs_noise}
\end{equation}

\subsubsection{FKE formulation}

Building upon the kinematically constrained tracking model defined above, we now map these specific system dynamics into the general Yau--Yau filtering framework to derive the exact governing FKE. Let the state variable be $x_t=[s(t),v(t)]^\top$ and the observation variable be $y_t=[R(t),\Theta(t)]^\top$. The objective is to derive the specific FKE that governs the spatiotemporal evolution of the probability density function $u(s,v,t)$. Assume diagonal process and observation covariance matrices
\begin{equation}
Q_t=\begin{bmatrix}Q_{11}&0\\0&Q_{22}\end{bmatrix},\qquad
S_t=\begin{bmatrix}S_{11}&0\\0&S_{22}\end{bmatrix}.
\label{eq:diag_cov}
\end{equation}

According to the DMZ equation theory in Section~\ref{subsec:yauyau}, the differential operator is defined by Eq.~\eqref{eq:operator}. Substituting the state variables $(x_1,x_2)=(s,v)$ and the drift functions $(f_1,f_2)=(v,f_v(s,t))$, and assuming the tangential acceleration is independent of velocity, the operator simplifies to
\begin{equation}
L(u)=\frac{1}{2}\left(Q_{11}\frac{\partial^2 u}{\partial s^2}
+Q_{22}\frac{\partial^2 u}{\partial v^2}\right)
-v\frac{\partial u}{\partial s}
-f_v\frac{\partial u}{\partial v}.
\label{eq:L_sv}
\end{equation}

By substituting the specific observation function vector $h=(h_R,h_\Theta)^\top$, the observation potential term becomes
\begin{equation}
\frac{1}{2}h^\top S^{-1}h
=\frac{1}{2}\left(\frac{h_R^2}{S_{11}}+\frac{h_\Theta^2}{S_{22}}\right).
\label{eq:obs_potential}
\end{equation}

Combining the operator and the observation potential term, the specific FKE for the fixed-trajectory target-tracking model can be explicitly written as
\begin{equation}
\frac{\partial u}{\partial t}
=
\frac{1}{2}\left(Q_{11}\frac{\partial^2 u}{\partial s^2}
+Q_{22}\frac{\partial^2 u}{\partial v^2}\right)
-v\frac{\partial u}{\partial s}
-f_v\frac{\partial u}{\partial v}
-\frac{1}{2}\left(\frac{h_R^2}{S_{11}}+\frac{h_\Theta^2}{S_{22}}\right)u.
\label{eq:fke_sv}
\end{equation}

Next, when the discrete observations arrive at time $\tau_i$, the initial condition for the subsequent FKE computation interval is  updated. Let the new observations be denoted by $R_i$ and $\Theta_i$. Expanding Eq.~\eqref{eq:yy_update} with the specific geometric functions yields the full initial-condition update equation:
\begin{equation}
u_i(s,v,\tau_{i-1})=
\exp\!\left[
\frac{R_i-R_{i-1}}{S_{11}}h_R(s)
+\frac{\Theta_i-\Theta_{i-1}}{S_{22}}h_\Theta(s)
\right]
u_{i-1}(s,v,\tau_{i-1}).
\label{eq:specific_update}
\end{equation}

Finally, after obtaining the updated unnormalized probability density at the current timestep, the state estimate is determined by calculating the expected value over the state space.

\subsection{Recentered-domain strategy}\label{subsec:recentered}

While Eq.~\eqref{eq:fke_sv} provides the mathematical foundation for target tracking, solving this FKE over a large global state space at every timestep can lead to high computational cost and limit real-time applicability.
Moving-window ideas have been widely used in recursive estimation, especially in moving-horizon estimation for constrained nonlinear systems \cite{ref28}. However, these methods shift a temporal estimation horizon rather than a spatial computational domain for FKE approximation. They therefore do not directly resolve the loss of spatial resolution that occurs when a fixed FKE domain contains large low-probability regions.
To address this bottleneck, we introduce a localized recentered-domain strategy.

The core of the strategy is to dynamically center a localized computational window on the state estimate obtained from the previous step. Fig.~\ref{fig:recentered} illustrates the fixed-size recentered-domain approach employed in RD-YYF for nonlinear filtering problems. This method aims to maintain computational efficiency and tracking accuracy by dynamically repositioning a compact computational domain around the most recent state estimate. The process consists of initial FKE computation, solution update, and domain transition to the next timestep.

\begin{inlinefigure}
\centering
\includegraphics[width=0.75\textwidth]{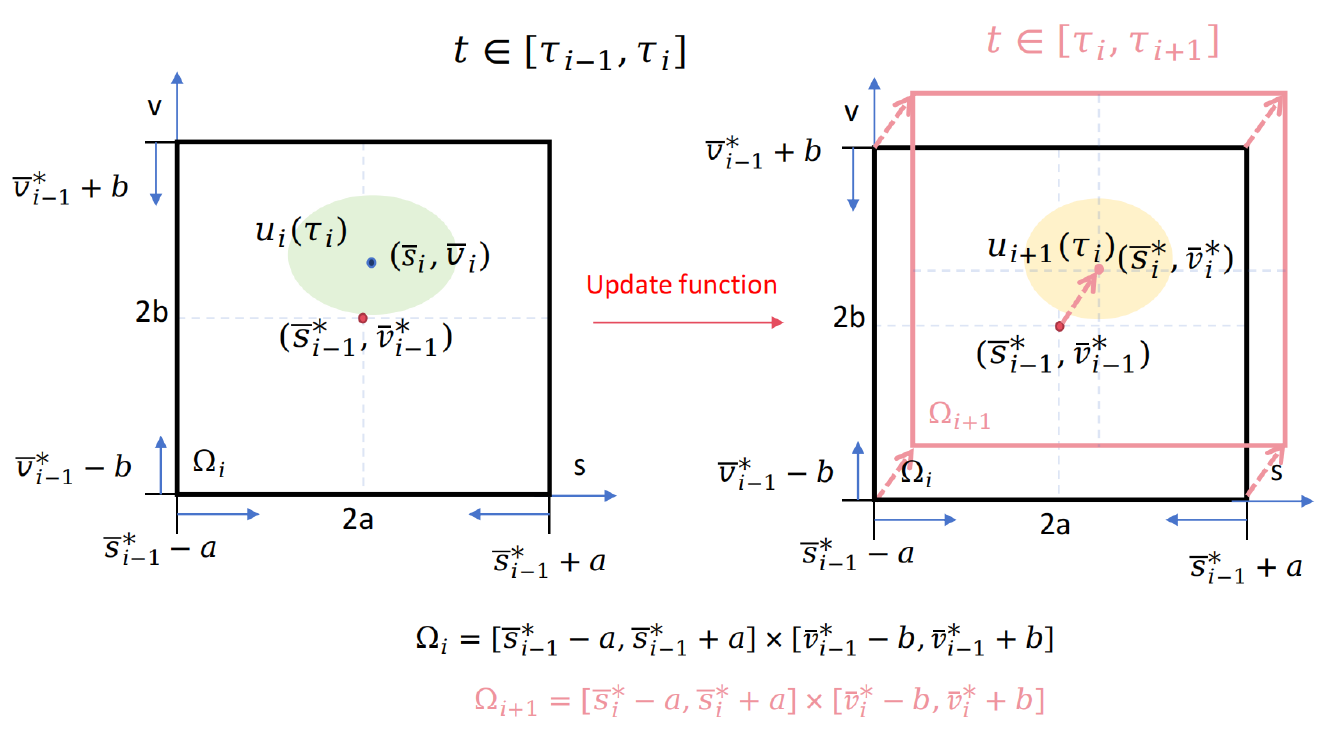}
\caption{Localized recentered-domain process, including the PINN solution, state update through the update function, and domain shifting for the next timestep.}
\label{fig:recentered}
\end{inlinefigure}

Specifically, the left panel of Fig.~\ref{fig:recentered} illustrates the FKE solving phase. The computational domain $\Omega_i$ spans $2a\times 2b$ in the arc-length and velocity dimensions. The domain is centered at the previous estimated state $(\bar{s}_{i-1}^*,\bar{v}_{i-1}^*)$, represented by the density center. The PINN computes the probability density $u_i$ within $\Omega_i$.

The middle panel depicts the observation update triggered by the new measurements $(R_i,\Theta_i)$. The density is  updated according to Eq.~\eqref{eq:specific_update}. Subsequently, the computational domain remains fixed for the current timestep, and the new estimated state $(\bar{s}_i^*,\bar{v}_i^*)$ is calculated using
\begin{equation}
\bar{s}_i^*=\frac{\int_{\Omega_i} s\,u_{i+1}(s,v,\tau_i)\,d\Omega_i}
{\int_{\Omega_i} u_{i+1}(s,v,\tau_i)\,d\Omega_i},\qquad
\bar{v}_i^*=\frac{\int_{\Omega_i} v\,u_{i+1}(s,v,\tau_i)\,d\Omega_i}
{\int_{\Omega_i} u_{i+1}(s,v,\tau_i)\,d\Omega_i}.
\label{eq:expectation}
\end{equation}

Finally, the right panel demonstrates the domain transition. The computational domain is shifted to $\Omega_{i+1}$ and aligned with the most recent state estimate $(\bar{s}_i^*,\bar{v}_i^*)$. To maintain temporal continuity for the next timestep, the initial condition on the shifted domain is obtained by transferring the updated density from the overlapping region of two consecutive computational domains:
\begin{equation}
u_{\mathrm{IC}}(s,v,\tau_i)=
\begin{cases}
u_{i+1}(s,v,\tau_i), & (s,v)\in \Omega_i\cap\Omega_{i+1},\\
0, & (s,v)\in \Omega_{i+1}\setminus\Omega_i.
\end{cases}
\label{eq:domain_transfer}
\end{equation}

By concentrating the computational domain around the latest state estimate, this localized approach reduces the need for a large global grid while maintaining numerical resolution near the dominant probability mass.

\subsection{Offline-online decoupled execution paradigm}\label{subsec:offline_online}

The proposed RD-YYF is implemented through an offline-online workflow. In the offline stage, PINN-generated FKE solution pairs are collected on recentered local domains, projected onto a PCA basis, and used to train a lightweight surrogate solver. In the online stage, the pretrained surrogate predicts the density evolution within the updated local window, followed by observation update and state estimation. The two stages are detailed in Sections~\ref{subsubsec:offline_stage} and~\ref{subsubsec:online_stage}.

\subsubsection{Offline training stage}\label{subsubsec:offline_stage}

The primary objective of the offline stage is to construct a lightweight FKE surrogate solver. This is achieved through a two-step process: generating a comprehensive database of high-fidelity FKE solutions using PINNs, followed by PCA-driven model reduction. The procedure is illustrated in Fig.~\ref{fig:offline}.

\begin{inlinefigure}
\centering
\includegraphics[width=0.8\textwidth]{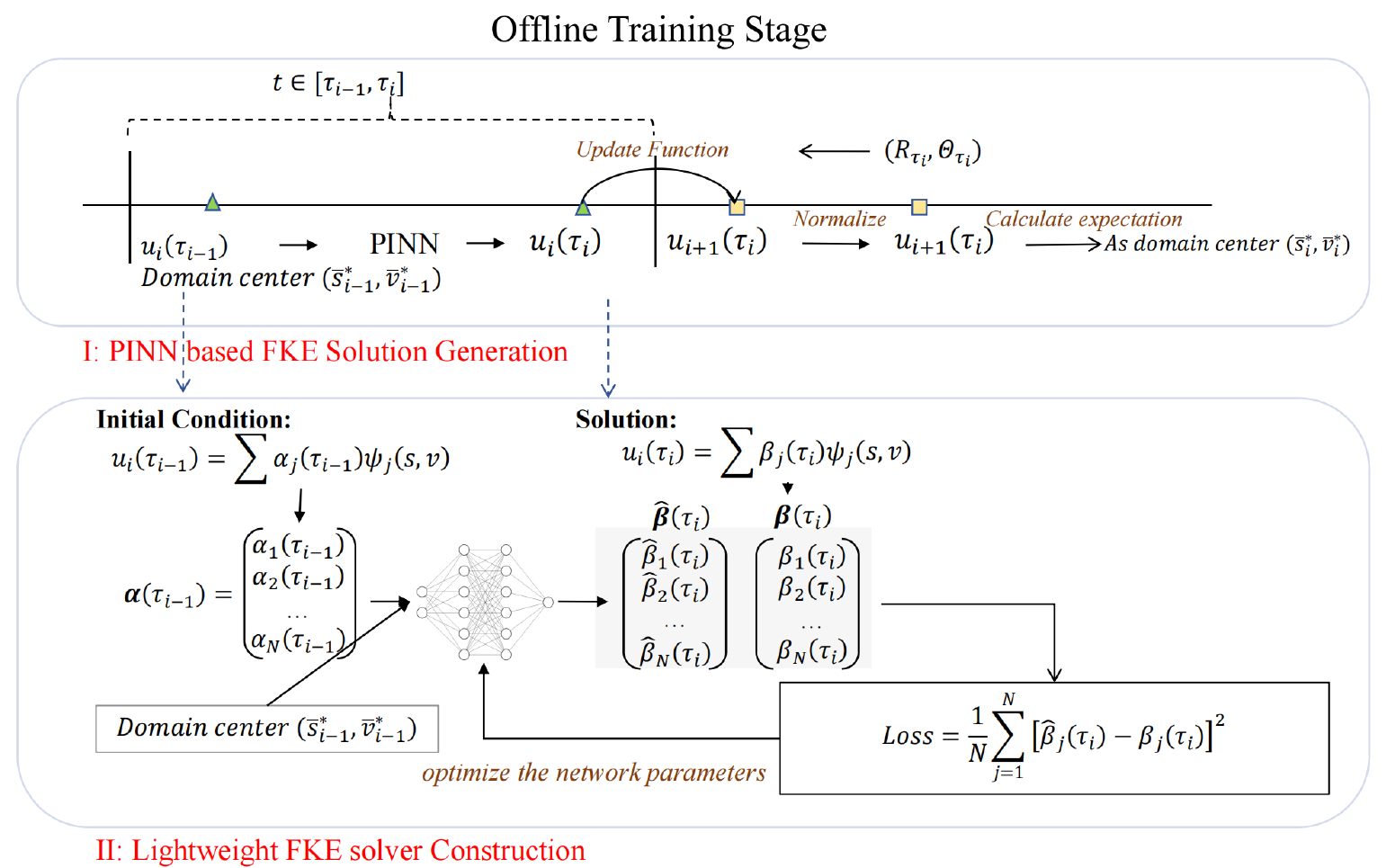}
\caption{The procedure of the offline training stage.}
\label{fig:offline}
\end{inlinefigure}

\textbf{I: PINN-based FKE solution generation.}
First, to build the foundational solution database, PINN is employed to solve the specific FKE sequentially over discrete time intervals. At each timestep within each localized computational domain $\Omega_i$, the network optimizes its parameters by minimizing the total loss
\begin{equation}
\mathcal{L}=\lambda_{\mathrm{FKE}}\mathcal{L}_{\mathrm{FKE}}
+\lambda_{\mathrm{IC}}\mathcal{L}_{\mathrm{IC}}
+\lambda_{\mathrm{BC}}\mathcal{L}_{\mathrm{BC}},
\label{eq:pinn_loss}
\end{equation}

Where $\lambda_{\mathrm{FKE}}$, $\lambda_{\mathrm{IC}}$, and $\lambda_{\mathrm{BC}}$ balance the relative importance of the FKE residual, initial-condition residual, and boundary-condition residual. The individual loss components are defined as
\begin{equation}
\begin{aligned}
\mathcal{L}_{\mathrm{FKE}}
&=\frac{1}{N_{\mathrm{FKE}}}
\sum_{j=1}^{N_{\mathrm{FKE}}}
\left|\mathcal{R}_{\mathrm{FKE}}\!\left(u_\theta(x^{j}_{\mathrm{FKE}},t^{j}_{\mathrm{FKE}})\right)\right|^2, (x^{j}_{\mathrm{FKE}},t^{j}_{\mathrm{FKE}})\sim\Omega_i\times[\tau_{i-1},\tau_i],\\
\mathcal{L}_{\mathrm{IC}}
&=\frac{1}{N_{\mathrm{IC}}}
\sum_{j=1}^{N_{\mathrm{IC}}}
\left|\mathcal{R}_{\mathrm{IC}}\!\left(u_\theta(x^{j}_{\mathrm{IC}},\tau_{i-1})\right)\right|^2,
 x^{j}_{\mathrm{IC}}\sim\Omega_i,\\
\mathcal{L}_{\mathrm{BC}}
&=\frac{1}{N_{\mathrm{BC}}}
\sum_{j=1}^{N_{\mathrm{BC}}}
\left|\mathcal{R}_{\mathrm{BC}}\!\left(u_\theta(x^{j}_{\mathrm{BC}},t^{j}_{\mathrm{BC}})\right)\right|^2,
 (x^{j}_{\mathrm{BC}},t^{j}_{\mathrm{BC}})\sim\partial\Omega_i\times[\tau_{i-1},\tau_i].
\end{aligned}
\label{eq:loss_components}
\end{equation}

Here, $\mathcal{L}_{\mathrm{FKE}}$, $\mathcal{L}_{\mathrm{IC}}$, and $\mathcal{L}_{\mathrm{BC}}$ denote the FKE, initial-condition, and boundary-condition residual losses, respectively. The values of $N_{\mathrm{FKE}}$, $N_{\mathrm{IC}}$, and $N_{\mathrm{BC}}$ are selected to ensure high accuracy.

During each training epoch, a new set of collocation points is generated to ensure uniform accuracy across the domain. The network parameters are optimized via gradient descent until the loss falls below the prescribed convergence thresholds or the number of epochs exceeds the maximum $E_{\max}$. Crucially, as illustrated in the upper panel of Fig.~\ref{fig:offline}, the offline stage simulates the step-by-step tracking workflow defined in Section~\ref{subsec:problem}. However, the ultimate goal here is data collection. For each interval $[\tau_{i-1},\tau_i]$, once the PINN successfully converges, the initial condition $u_i(\tau_{i-1})$ and the corresponding end-time solution $u_i(\tau_i)$ are extracted and stored. This iterative process constructs the high-dimensional spatiotemporal dataset required for the subsequent dimensionality reduction.

The PINN-based database generation procedures are summarized in Algorithms~\ref{alg:offline_pinn}.

\begin{algorithm}[H]
\caption{PINN-based FKE solution generation}
\label{alg:offline_pinn}
\begin{algorithmic}[1]
\State \textbf{Input:} Number of timesteps $N_T$, interval $\Delta t=[\tau_{i-1},\tau_i]$, convergence thresholds $\epsilon=(\epsilon_{\mathrm{FKE}},\epsilon_{\mathrm{IC}},\epsilon_{\mathrm{BC}})$, maximum epochs $E_{\max}$, and observations $(R_{\tau_i},\Theta_{\tau_i})$.
\State Randomly initialize PINN parameters $\theta$, set $\Omega_0$ centered at $(s_0,v_0)$, and set the initial density $u_1(s,v,\tau_0)$.
\For{$i=1$ to $N_T$}
    \For{epoch $=1$ to $E_{\max}$}
        \State Solve Eq.~\eqref{eq:fke_sv} on $\Omega_{i-1}$ with initial condition $u_i(s,v,\tau_{i-1})$.
        \State Compute
        $\displaystyle
        \mathcal{L}$
        \State Update the network parameters $\theta$ by gradient descent.
        \If{$\mathcal{L}<\epsilon$}
            \State break.
        \EndIf
    \EndFor
    \State Evaluate the terminal FKE solution $u_i(s,v,\tau_i)$.
    \State Apply the observation update in Eq.~\eqref{eq:specific_update} and normalize $u_{i+1}(s,v,\tau_i)$.
    \State Compute $(\bar{s}_i^*,\bar{v}_i^*)$ by Eq.~\eqref{eq:expectation}, recenter $\Omega_{i+1}$, and transfer the initial condition by Eq.~\eqref{eq:domain_transfer}.
\EndFor
\State \textbf{Output:} High-dimensional FKE solution-pair dataset
\Statex \hspace{\algorithmicindent}$\displaystyle
\mathcal{S}=\{u_i(s,v,\tau_{i-1}),u_i(s,v,\tau_i),(\bar{s}_{i-1}^*,\bar{v}_{i-1}^*)\}_{i=1}^{N_T}.$
\end{algorithmic}
\end{algorithm}

\textbf{II: PCA-driven model reduction.}
Second, PCA is employed to extract dominant modes from the FKE solution snapshots generated in the previous stage. These PCs provide a low-dimensional basis that preserves the main variance of the density evolution while reducing the representation dimension. The corresponding basis functions $\psi(s,v)$ form an orthonormal set for representing the FKE solution snapshots. For each time interval, the initial conditions and terminal solutions can be expressed as linear combinations of the basis functions:
\begin{equation}
u_i(\tau_{i-1})=\alpha^\top(\tau_i)\psi(s,v),\qquad
u_i(\tau_i)=\beta^\top(\tau_i)\psi(s,v),
\label{eq:pca_rep}
\end{equation}

Where $\alpha$ and $\beta$ are coefficient vectors representing the projection weights onto each PC. These extracted coefficients are subsequently used to train the lightweight FKE solver. To account for the dynamically updated window, the solver is designed to learn the nonlinear mapping from the concatenated inputs, comprising the initial-condition coefficients and the corresponding absolute domain center, to the target solution coefficients. This training process minimizes the loss between the true $\beta(\tau_i)$ and predicted $\widehat{\beta}(\tau_i)$:
\begin{equation}
\mathcal{L}_{\mathrm{sur}}=\left\|\beta(\tau_i)-\widehat{\beta}(\tau_i)\right\|_2^2.
\label{eq:sur_loss}
\end{equation}

The specific training procedure for this lightweight surrogate network is detailed in Algorithms~\ref{alg:surrogate}.

\begin{algorithm}[H]
\caption{Lightweight FKE solver construction}
\label{alg:surrogate}
\begin{algorithmic}[1]
\State \textbf{Input:} Generated dataset $\mathcal{S}$, number of retained PCs $N$, and maximum epochs $E_{\max}$.
\State Randomly initialize surrogate parameters $\eta$.
\State Perform PCA on $\mathcal{S}$ and retain the leading components $\{\psi_j(s,v)\}_{j=1}^{N}$.
\State Project each initial and terminal FKE solution onto the PCA basis:
\Statex \hspace{\algorithmicindent}$\displaystyle
u_i(s,v,\tau_{i-1})\mapsto
\alpha(\tau_{i-1})=[\alpha_1(\tau_{i-1}),\ldots,\alpha_N(\tau_{i-1})]^\top,$
\Statex \hspace{\algorithmicindent}$\displaystyle
u_i(s,v,\tau_i)\mapsto
\beta(\tau_i)=[\beta_1(\tau_i),\ldots,\beta_N(\tau_i)]^\top.$
\For{epoch $=1$ to $E_{\max}$}
    \State Predict terminal coefficients using
    \Statex \hspace{\algorithmicindent}$\displaystyle
    \widehat{\beta}(\tau_i)=
    F_{\eta}\!\left(\alpha(\tau_{i-1}),(\bar{s}_{i-1}^*,\bar{v}_{i-1}^*)\right).$
    \State Minimize
    \Statex \hspace{\algorithmicindent}$\displaystyle
    \mathcal{L}_{\mathrm{sur}}
    =\|\beta(\tau_i)-\widehat{\beta}(\tau_i)\|_2^2,$
    \State and optimize $\eta$.
\EndFor
\State \textbf{Output:} Trained lightweight FKE solver.
\end{algorithmic}
\end{algorithm}

Algorithms~\ref{alg:offline_pinn} and~\ref{alg:surrogate} involve iterative PINN optimization, PCA construction, and surrogate training, and are therefore performed in the offline stage. After this offline construction is completed, the pretrained lightweight FKE solver is deployed in the online stage for per-timestep density prediction and state estimation.

\subsubsection{Online execution stage}\label{subsubsec:online_stage}

In the online execution stage, the lightweight FKE solver is deployed for per-timestep density prediction and state estimation. As illustrated in Fig.~\ref{fig:online}, the sequential state evolution is propagated through forward passes of the pretrained surrogate network.

\begin{inlinefigure}
\centering
\includegraphics[width=\textwidth]{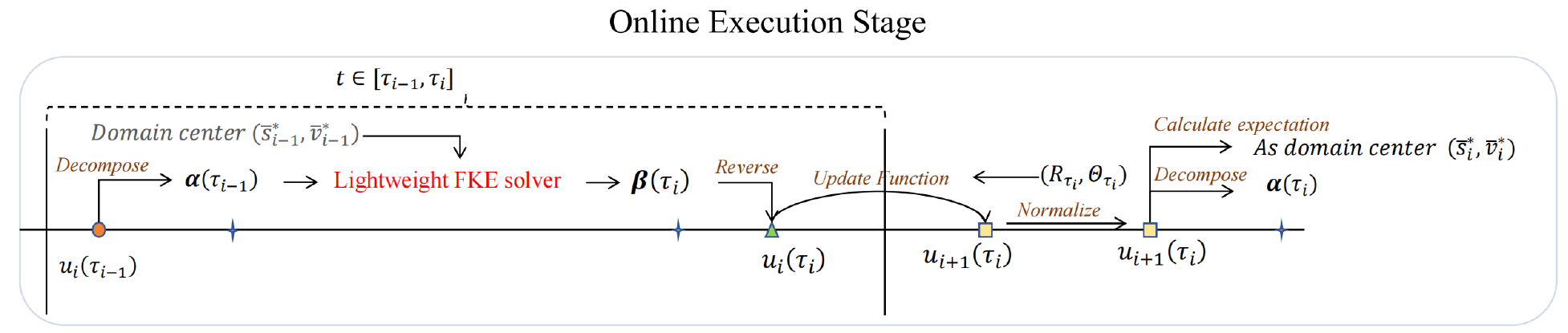}
\caption{The process of the online execution stage.}
\label{fig:online}
\end{inlinefigure}

For each incoming timestep, the initial probability density $u_i(\tau_{i-1})$ within the localized window is first mapped to the standardized local grid and projected onto the orthogonal PCA basis to extract the low-dimensional initial-condition coefficients $\alpha(\tau_{i-1})$. These coefficients, concatenated with the current domain center coordinates $(\bar{s}_{i-1}^*,\bar{v}_{i-1}^*)$, are fed directly into the surrogate solver to predict the end-time solution coefficients $\beta(\tau_i)$. The end solution density $u_i(\tau_i)$ is then rapidly reconstructed via a linear combination of the PCs and mapped back to the domain centered at $(\bar{s}_{i-1}^*,\bar{v}_{i-1}^*)$. Upon receiving the real-time measurement $(R_{\tau_i},\Theta_{\tau_i})$, the update function is applied to yield the posterior density. Finally, the new state expectation is computed, and the localized computational window is dynamically recentered for the subsequent timestep.
The online RD-YYF execution procedure is summarized in Algorithm~\ref{alg:online}.

\begin{algorithm}[H]
\caption{Online RD-YYF execution}
\label{alg:online}
\begin{algorithmic}[1]
\State \textbf{Input:} PCA basis $\{\psi_j(s,v)\}_{j=1}^{N}$, lightweight FKE solver $F_{\eta}$, and observations $(R_{\tau_i},\Theta_{\tau_i})$.
\State Initialize $\Omega_0$ centered at $(s_0,v_0)$ and set the initial density $u_1(s,v,\tau_0)$.
\For{$i=1$ to $N_T$}
    \State Map the initial condition to the standardized local grid and project it onto the PCA basis:
    \Statex \hspace{\algorithmicindent}$\displaystyle
    \alpha_j(\tau_{i-1})
    =\langle u_i(s,v,\tau_{i-1}),\psi_j(s,v)\rangle,\quad j=1,\ldots,N.$
    \State Predict the terminal coefficients:
    \Statex \hspace{\algorithmicindent}$\displaystyle
    \widehat{\beta}(\tau_i)=
    F_{\eta}\!\left(\alpha(\tau_{i-1}),(\bar{s}_{i-1}^*,\bar{v}_{i-1}^*)\right).$
    \State Reconstruct the terminal density:
    \Statex \hspace{\algorithmicindent}$\displaystyle
    \widehat{u}_i(s,v,\tau_i)=\widehat{\beta}^{\top}(\tau_i)\psi(s,v).$
    \State Apply the observation update in Eq.~\eqref{eq:specific_update}.
    \State Normalize $u_{i+1}(s,v,\tau_i)$, compute $(\bar{s}_i^*,\bar{v}_i^*)$ by Eq.~\eqref{eq:expectation}, and recenter $\Omega_{i+1}$.
\EndFor
\end{algorithmic}
\end{algorithm}

This offline-online RD-YYF framework combines PINN-based PDE solution generation with PCA-based feature extraction and dimensionality reduction. The offline training stages prepare the lightweight solver, enabling rapid online computation of FKE solutions for nonlinear filtering problems.

\section{Simulations}\label{sec:simulations}

This section evaluates the tracking accuracy and online computational behavior of the proposed RD-YYF algorithm through numerical simulations. To provide a comparative analysis, EKF and PF are selected as baseline methods because they are widely used in nonlinear target tracking and represent two standard filtering paradigms: linearization-based Gaussian approximation and Monte Carlo-based Bayesian estimation, respectively \cite{ref5,ref7,ref29}. These baselines provide a meaningful reference for evaluating both tracking accuracy and computational cost. All experimental scenarios are formulated under the kinematically constrained continuous-time tracking framework established in Section~\ref{subsec:problem}.

\subsection{Implementation details and experimental setup}\label{subsec:setup}

This subsection details the specific neural network architectures and hyperparameter configurations employed during the two-stage offline training procedure, alongside the computational environment and evaluation metrics used for real-time tracking.

\subsubsection{Network architectures and configurations}

In both the offline and online stages, the localized computational window uses the same fixed size. Unless otherwise specified, the half-widths of the local domain are set to 4 in all simulations. This setting was chosen based on preliminary numerical tests to cover the dominant support of the posterior density in the tested scenarios while keeping the FKE computation localized.

To separate offline surrogate construction from online performance evaluation, independent random seeds were used for the two stages. For each example, 14 Monte Carlo trajectories were generated in the offline stage for PINN snapshot generation, PCA construction, and surrogate training, whereas 20 additional Monte Carlo trajectories with non-overlapping seeds were used for online MSE evaluation.

In offline training stage I, the PINN is constructed as a fully connected neural network utilizing the hyperbolic tangent activation function, which effectively captures nonlinear dynamics. The network parameters are optimized using the Adam optimizer \cite{ref30} to balance convergence speed and numerical stability.

As detailed in Table~\ref{tab:pinn}, distinct convergence thresholds and penalty weights are employed to account for the different loss components. This dynamic balancing of physical constraints prevents overfitting and limits excessive computational cost, while ensuring that the generated PDE solutions meet the precision required for the subsequent dimensionality reduction.

\begin{inlinetable}
\caption{PINN configuration for offline training.}
\label{tab:pinn}
\centering
\begin{tabular}{@{}ll@{}}
\toprule
Parameter & Value\\
\midrule
Number of hidden layers & 3\\
Nodes per layer & 100\\
Learning rate & 0.008\\
Max epochs per timestep & 50000\\
Convergence thresholds & $\epsilon_{\mathrm{FKE}}=2\times10^{-4}$, $\epsilon_{\mathrm{IC}}=5\times10^{-6}$, $\epsilon_{\mathrm{BC}}=1\times10^{-5}$\\
\bottomrule
\end{tabular}
\end{inlinetable}

In offline training stage II, the number of retained PCs is determined according to two criteria: the cumulative explained variance of the FKE solution snapshots and the validation performance of the lightweight solver.The retained dimension is selected to achieve a compact representation while maintaining stable validation accuracy.

Since the solution at each timestep differs only slightly from the initial condition, particularly when the time interval is small, the network learns the difference between the PC coefficients rather than learning the full mapping directly. Following the residual-learning idea \cite{ref31}, an adapted residual network structure is adopted as the approximate FKE solver. The residual network consists of several residual blocks, where each residual block contains a fully connected network with a skip connection. The skip connection adds the block input directly to its output, enhancing gradient flow during backpropagation and stabilizing offline training.

In this work, the residual network is configured with three residual blocks, where each block contains two fully connected layers with 64 nodes per layer. Training is performed using the Adam optimizer with a learning rate of 0.001 for 5000 epochs to ensure convergence.

\subsubsection{Computational environment and performance metrics}

All numerical simulations were conducted on an NVIDIA RTX3070 GPU and a platform with 16 Intel Core i5-13600KF CPU cores at 3.50 GHz. The neural networks were implemented using PyTorch.

The reported CPU time refers to the average per-timestep online filtering cost after offline training. For RD-YYF, it includes PCA projection, surrogate prediction, density reconstruction, observation update, normalization, and state-expectation computation, but excludes offline PINN snapshot generation, PCA construction, and surrogate training. EKF and PF are timed using the same per-timestep online protocol.

The tracking performance is evaluated using the mean squared error (MSE) for each state variable:
\begin{equation}
\mathrm{MSE}_{x_d}=\frac{1}{N_T}\sum_{i=1}^{N_T}\left(x_d(\tau_i)-\widehat{x}_d(\tau_i)\right)^2,
\label{eq:mse}
\end{equation}
where $x_d(\tau_i)$ is the true value of the $d$th state component at time $\tau_i$, and $\widehat{x}_d(\tau_i)$ is the corresponding online estimate.

\subsection{Target-tracking performance}\label{subsec:performance}

To evaluate the tracking performance of the proposed RD-YYF, we consider two representative target-tracking scenarios within the continuous-time, kinematically constrained framework established in Section~\ref{subsec:problem}. The two examples differ mainly in the predefined geometric trajectory of the target, allowing the proposed method to be assessed under both simple and curved motion constraints. 

\subsubsection{Example 1: straight-line trajectory}

In this example, we consider a target-tracking scenario along a straight-line trajectory, which serves as a linear geometric constraint within the established framework. This model acts as a mathematical abstraction of a constrained target monitored via nonlinear range and bearing measurements, a scenario commonly encountered in rail-based transportation systems or guided motion platforms.

Assuming the straight-line path aligns with the Cartesian $x$-axis, the physical arc-length state directly corresponds to the one-dimensional position. An observer is situated at a fixed off-path location $(0,1)$. The objective is to accurately track the coupled kinematic states, namely the target position $s(t)$ along the path and its velocity $v(t)$, using continuous noisy range and bearing observations. The system dynamics and observation processes are described by
\begin{equation}
\left\{
\begin{aligned}
ds(t)&=v(t)\,dt+d\nu_s(t),\\
dv(t)&=-\left(\frac{\pi}{2}\right)^2s(t)\,dt+d\nu_v(t),\\
dR(t)&=\sqrt{s^2(t)+1}\,dt+d w_R(t),\\
d\Theta(t)&=\arctan(s(t))\,dt+d w_\Theta(t),
\end{aligned}
\right.
\label{eq:ex1_model}
\end{equation}

Where $\nu_t$ and $w_t$ are Brownian motion processes with $E[d\nu_t d\nu_t^\top]=0.1I_2\,dt$ and $E[d w_t d w_t^\top]=0.1I_2\,dt$, and $I_2$ is the two-dimensional identity matrix. The true initial physical states of the target are set to $s_0=0$ m and $v_0=1$ m/s. The continuous tracking experiment is conducted over a total simulation duration of $T=10$ s, with a discrete temporal resolution of $\Delta t=0.01$.

The corresponding FKE is
\begin{equation}
\frac{\partial u}{\partial t}
=0.05\left(\frac{\partial^2 u}{\partial s^2}
+\frac{\partial^2 u}{\partial v^2}\right)
-v\frac{\partial u}{\partial s}
+\left(\frac{\pi}{2}\right)^2s\frac{\partial u}{\partial v}
-5\left(s^2+1+\arctan^2(s)\right)u.
\label{eq:ex1_fke}
\end{equation}

We assume the initial distribution for the first timestep is $u_0(s,v,0)=\exp[-2(s^2+(v-1)^2)]$. The localized computational domain for solving the FKE is configured with half-lengths $a=4$ and $b=4$. Centered around the initial true state $(s_0,v_0)$, the starting domain is precisely bounded as $\Omega_0=[-4,4]\times[-3,5]$. At each subsequent timestep, this compact window dynamically shifts to align with the latest state estimate.

To construct the offline database for PCA reduction, 14 Monte Carlo trajectories are simulated using the PINN-based FKE solver. The resulting high-fidelity FKE solution pairs are used for PCA-based model reduction, with average offline MSE values of MSE-$s=0.125$ and MSE-$v=0.315$.

Fig.~\ref{fig:epoch} illustrates the training behavior of the offline PINN solver over the 10-second tracking horizon for a representative simulation run. Each plotted loss value corresponds to the final converged loss at the terminal epoch of one discrete timestep. The required number of training epochs is high at the first timestep, but decreases rapidly and remains at a lower level in the subsequent steps. This trend is attributed to the continuous evolution of the probability density, which allows the optimized network weights from the previous timestep to provide an effective initialization for the current timestep. The final converged FKE, IC, and BC losses also remain consistently low across the time horizon. This behavior indicates stable PINN convergence during the sequential FKE computations. It also supports the effectiveness of the localized window strategy in maintaining the density approximation over successive timesteps.

\begin{inlinefigure}
\centering
\includegraphics[width=0.75\textwidth]{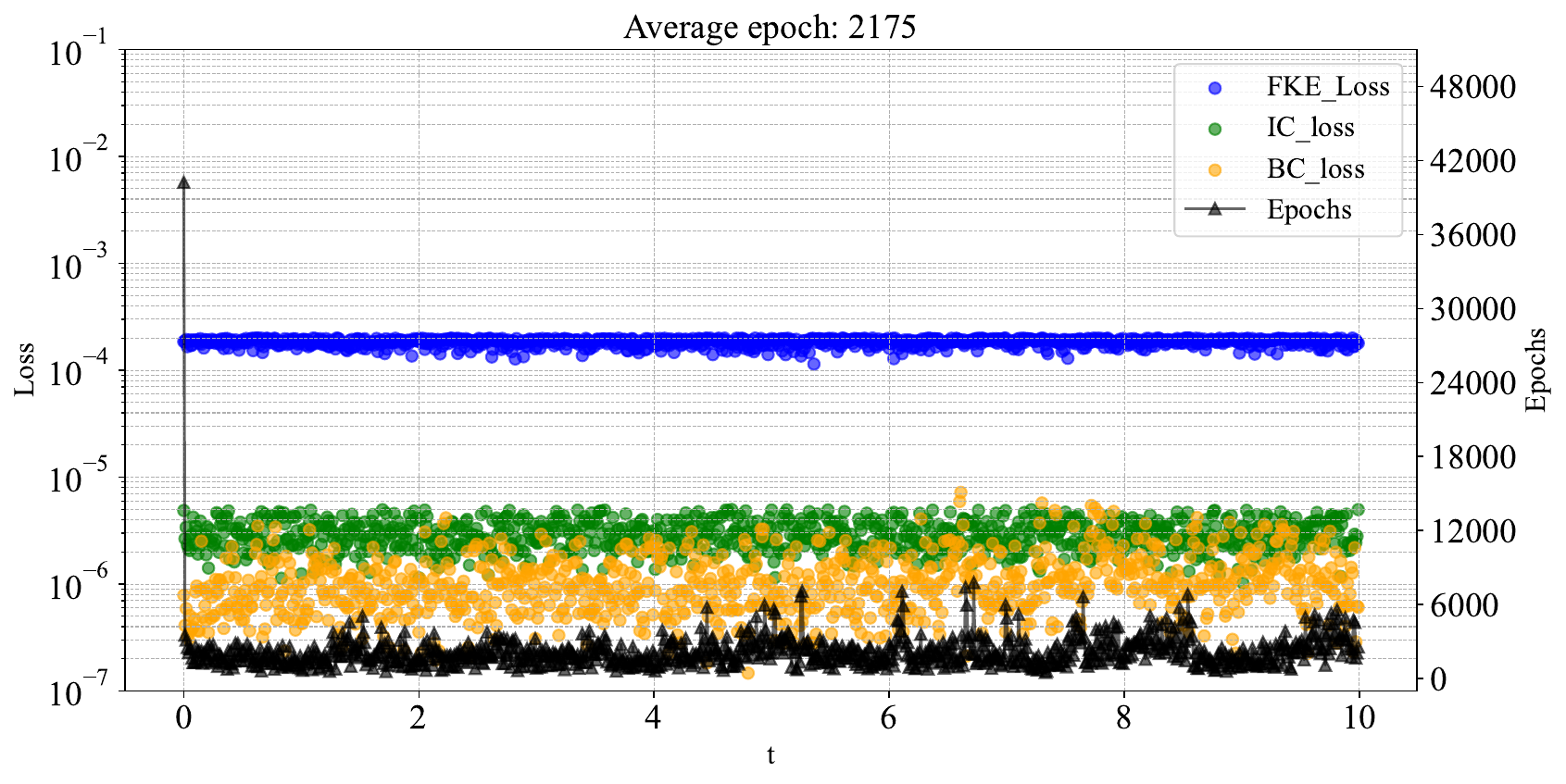}
\caption{Convergence behavior of the offline PINN solver across sequential timesteps for a representative simulation run. The scatter points denote the final converged values of the FKE, IC, and BC losses at each timestep, together with the total number of epochs required to reach convergence.}
\label{fig:epoch}
\end{inlinefigure}

During offline training stage II, the top 10 PCs are retained, capturing 96.97\% of the total variance (Fig.~\ref{fig:pca}).

\begin{inlinefigure}
\centering
\includegraphics[width=0.8\textwidth]{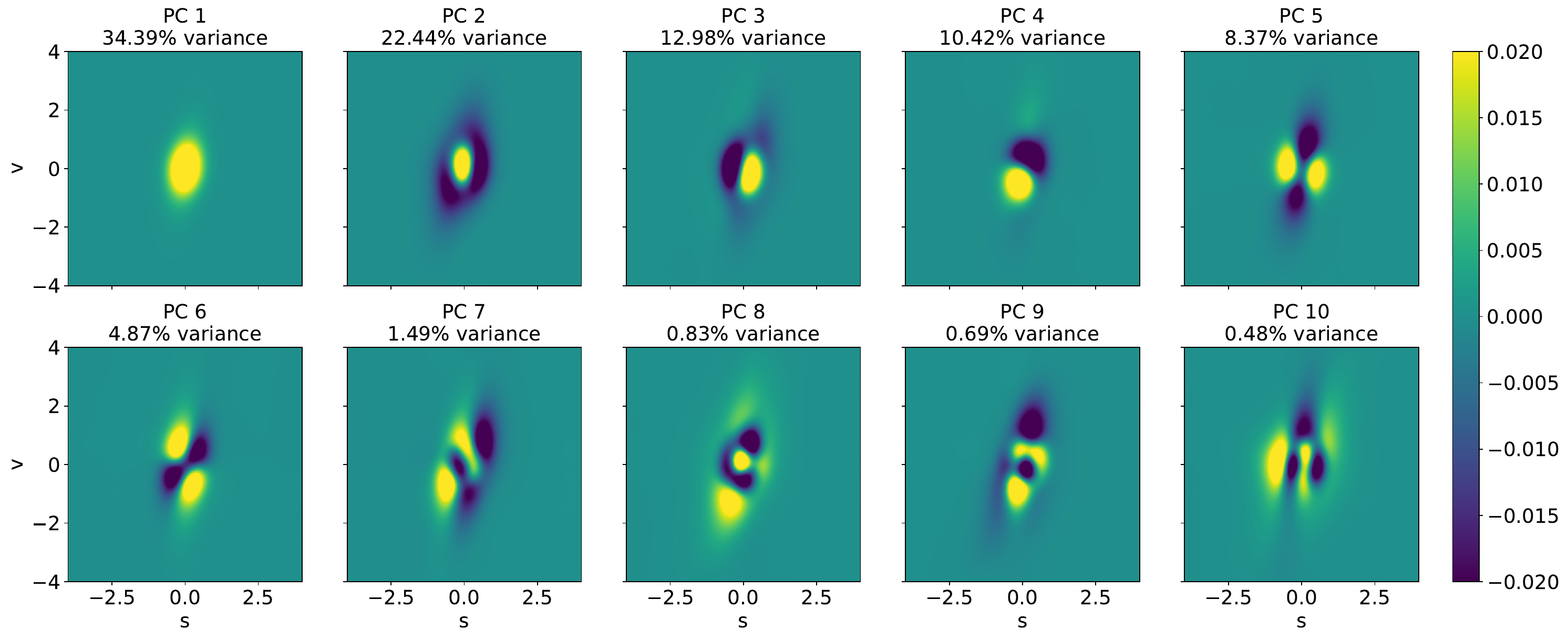}
\caption{The top 10 PCs and their explained variance ratios in Example 1.}
\label{fig:pca}
\end{inlinefigure}

The lightweight FKE solver is trained on these 10 component coefficients using 10,640 timesteps for training and 2,660 timesteps for testing within the offline database. This split is used to monitor the surrogate approximation capability, whereas final tracking performance is evaluated through the independent online filtering experiments reported in Table~\ref{tab:ex1}. The corresponding training and testing loss curves are shown in Fig.~\ref{fig:loss}.

\begin{inlinefigure}
\centering
\includegraphics[width=0.4\textwidth]{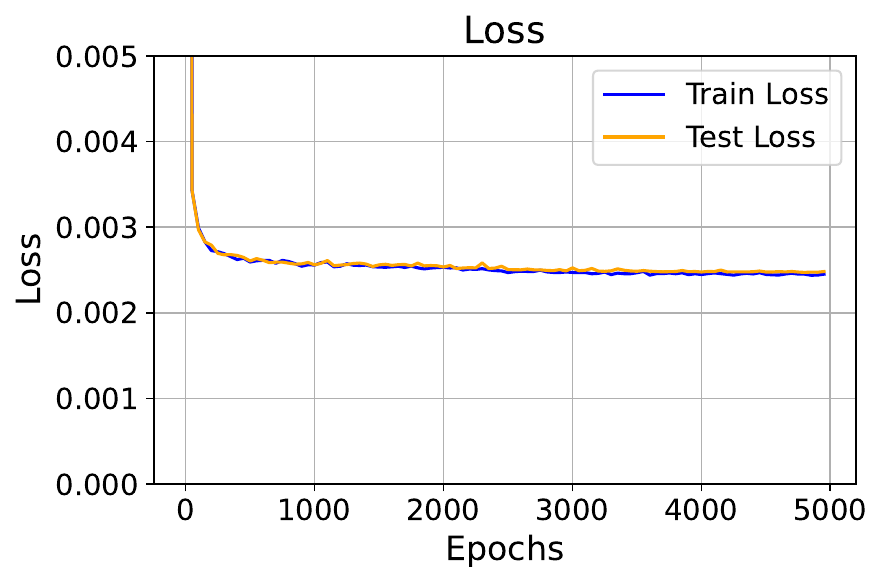}
\caption{Training and testing loss convergence of the lightweight FKE solver with 10 PCs.}
\label{fig:loss}
\end{inlinefigure}

In the online execution stage, the pretrained lightweight solver processes real-time measurement signals to predict the target state variables $s(t)$ and $v(t)$. To quantitatively evaluate the performance of the proposed RD-YYF, it is benchmarked against EKF and PF. For PF, configurations with both 1,000 and 10,000 particles are implemented to explicitly assess the impact of sample size on estimation accuracy.

Table~\ref{tab:ex1} reports the mean MSEs over 20 online Monte Carlo trajectories, together with the average per-timestep online CPU time. Increasing the number of PF particles from 1,000 to 10,000 does not lead to a clear improvement in this example. In contrast, RD-YYF achieves lower mean MSEs for both position and velocity, indicating better tracking accuracy under the tested setting. 

\begin{inlinetable}
\caption{Mean performance over 20 online Monte Carlo trajectories for EKF, PF with 1,000 and 10,000 particles, and RD-YYF in Example 1.}
\label{tab:ex1}
\centering
\begin{tabular}{@{}lccc@{}}
\toprule
Method & MSE-$s$ & MSE-$v$ & CPU time\\
\midrule
EKF & 0.178 & 0.512 & 0.066 ms\\
PF (1000) & 0.149 & 0.434 & 1.169 ms\\
PF (10000) & 0.147 & 0.433 & 1.470 ms\\
\textbf{RD-YYF} & \textbf{0.139} & \textbf{0.380} & 1.446 ms\\
\bottomrule
\end{tabular}
\end{inlinetable}

Fig.~\ref{fig:ex1} presents the estimated state trajectories from one representative online trial among the 20 Monte Carlo trajectories. The estimated position $s(t)$ and velocity $v(t)$ generated by EKF, PF, and RD-YYF are plotted alongside the ground truth. RD-YYF closely tracks the true state trajectories across the time horizon and provides smooth estimates in this representative case.

\begin{inlinefigure}
\centering
\includegraphics[width=0.9\textwidth]{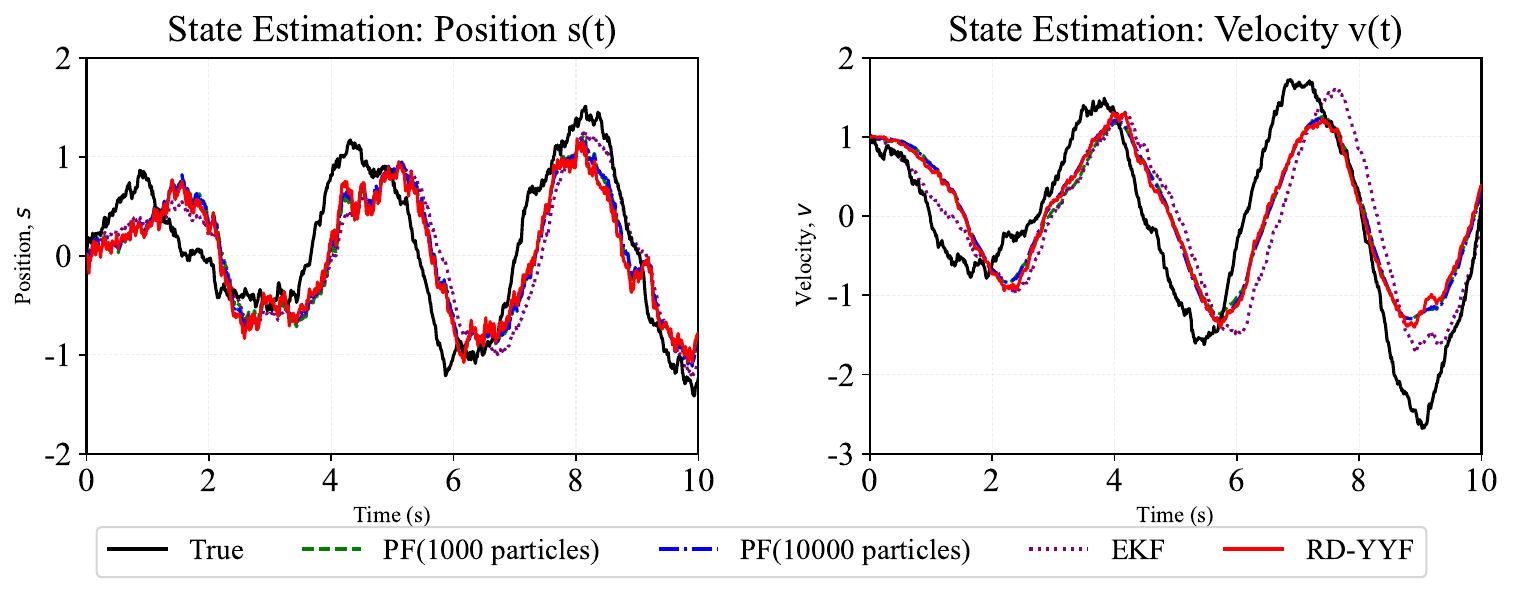}
\caption{Estimated state trajectories of $s(t)$ and $v(t)$ obtained by EKF, PF, and RD-YYF in one representative online Monte Carlo trajectory, compared with the ground truth.}
\label{fig:ex1}
\end{inlinefigure}

\FloatBarrier

\subsubsection{Example 2: curved trajectory}

To evaluate the proposed filter under more stringent conditions, the second scenario examines a target navigating a known curved path, characteristic of autonomous vehicles or mobile robots maneuvering along complex winding routes. Specifically, the target is constrained to move along a unit circular arc, whose parametric representation is
\begin{equation}
x=\phi_x(s)=\cos(s),\qquad y=\phi_y(s)=\sin(s).
\label{eq:circle}
\end{equation}

Where $s$ denotes the arc-length parameter along the circular trajectory. An external observer is fixed at Cartesian coordinates $(0,-2)$. The target dynamics are described in terms of the evolution of $s(t)$ and the corresponding tangential velocity $v(t)$, leading to
\begin{equation}
\left\{
\begin{aligned}
ds(t)&=v(t)\,dt+d\nu_s(t),\\
dv(t)&=-\pi^2s(t)\,dt+d\nu_v(t),\\
dR(t)&=\sqrt{\cos^2(s(t))+(\sin(s(t))+2)^2}\,dt+d w_R(t),\\
d\Theta(t)&=\arctan\frac{\cos(s(t))}{\sin(s(t))+2}\,dt+d w_\Theta(t).
\end{aligned}
\right.
\label{eq:ex2_model}
\end{equation}

Here, $\nu_t$ and $w_t$ are Brownian motion processes with process and observation noise covariances $0.1I_2\,dt$. The true initial position and velocity of the target are $s_0=0$ m and $v_0=0.5$ m/s. The simulation spans a total duration of $T=10$ s with time interval $\Delta t=0.01$.

The corresponding FKE is
\begin{equation}
\begin{aligned}
\frac{\partial u}{\partial t}
&=0.05\left(\frac{\partial^2 u}{\partial s^2}
+\frac{\partial^2 u}{\partial v^2}\right)
-v\frac{\partial u}{\partial s}
+\pi^2s\frac{\partial u}{\partial v} -5\left(\cos^2(s)+(\sin(s)+2)^2
+\arctan^2\frac{\cos(s)}{\sin(s)+2}\right)u.
\end{aligned}
\label{eq:ex2_fke}
\end{equation}

We assume the initial distribution for the first timestep is $u_0(s,v,0)=\exp[-2(s^2+(v-0.5)^2)]$. Following the identical offline training protocol established in Example 1, the localized computational window maintains its dimensions. Centered at the new initial state $(s_0,v_0)$, the initial domain is defined as $\Omega_0=[-4,4]\times[-3.5,5.5]$.
The average offline MSEs are MSE-$s=0.158$ and MSE-$v=0.217$. For PCA dimensionality reduction in Stage II, the top 15 PCs are retained, capturing 98.89\% of the total variance. The lightweight surrogate solver is subsequently trained on these 15 components using 10,640 timesteps for training and 2,660 timesteps for testing within the offline database, with a batch size of 512. 

Table~\ref{tab:ex2} summarizes the mean tracking performance over 20 online Monte Carlo trajectories. RD-YYF achieves the lowest mean MSEs for both position and velocity among the compared methods. Compared with PF using 10,000 particles, the gain in position estimation is marginal, whereas the gain in velocity estimation is more evident. For online computation, RD-YYF requires more per-timestep CPU time than EKF and PF with 1,000 particles, but less time than PF with 10,000 particles. These results indicate that RD-YYF provides lower average tracking errors under the tested settings, with a per-timestep online cost between the two PF configurations.

\begin{inlinetable}
\caption{Mean performance over 20 online Monte Carlo trajectories for EKF, PF with 1,000 and 10,000 particles, and RD-YYF in Example 2.}
\label{tab:ex2}
\centering
\begin{tabular}{@{}lccc@{}}
\toprule
Method & MSE-$s$ & MSE-$v$ & CPU time\\
\midrule
EKF & 0.162 & 0.248 & 0.036 ms\\
PF (1000) & 0.146 & 0.240 & 1.417 ms\\
PF (10000) & 0.143 & 0.231 & 2.100 ms\\
\textbf{RD-YYF} & \textbf{0.141} & \textbf{0.221} & 1.749 ms\\
\bottomrule
\end{tabular}
\end{inlinetable}

Fig.~\ref{fig:ex2} shows one representative online realization from the 20 Monte Carlo trajectories in Example 2. The estimated trajectories of $s(t)$ and $v(t)$ produced by EKF, PF, and RD-YYF are compared with the ground truth. In this representative curved-path case, RD-YYF maintains close agreement with the true trajectories under nonlinear dynamics and noisy observations.

\begin{inlinefigure}
\centering
\includegraphics[width=0.8\textwidth]{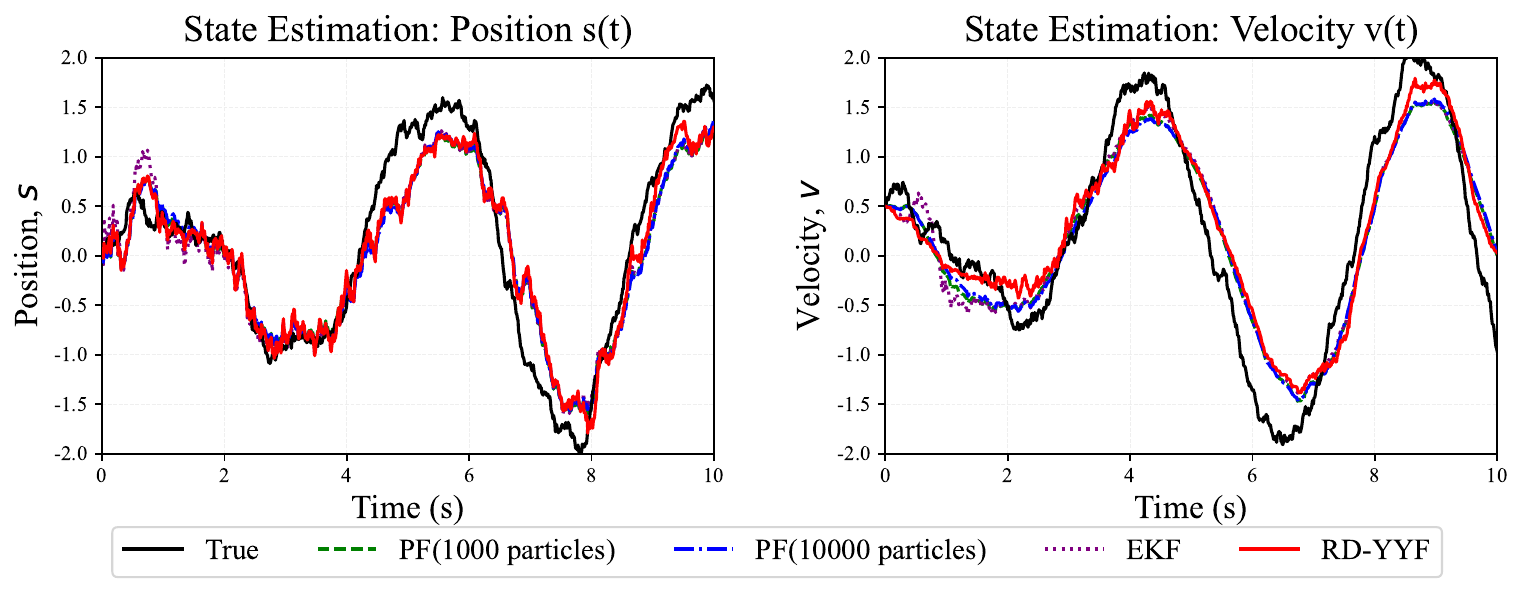}
\caption{Estimated state trajectories of $s(t)$ and $v(t)$ obtained by EKF, PF, and RD-YYF in one representative online Monte Carlo trajectory, compared with the ground truth.}
\label{fig:ex2}
\end{inlinefigure}

Fig.~\ref{fig:cartesian} compares the recovered Cartesian trajectories obtained by EKF, PF with 1,000 and 10,000 particles, and RD-YYF with the ground truth. The results show that RD-YYF remains effective after transformation from the arc-length representation to the original two-dimensional tracking space.

\begin{inlinefigure}
\centering
\includegraphics[width=0.8\textwidth]{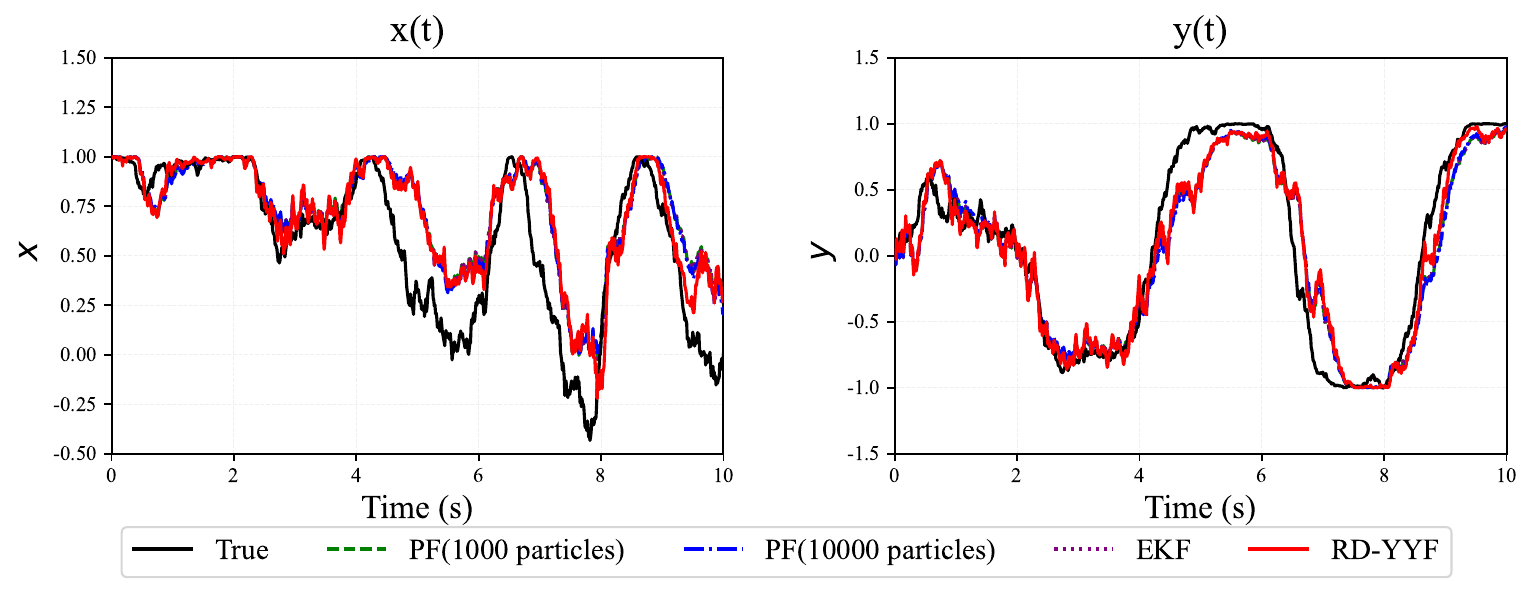}
\caption{Cartesian coordinate recovery results in Example 2.}
\label{fig:cartesian}
\end{inlinefigure}

\FloatBarrier

\subsection{Ablation study: effect of the recentered-domain strategy}\label{subsec:ablation}

To evaluate the necessity of the proposed recentered-domain strategy, we conduct an ablation study based on Example 2, which represents a more challenging geometrically constrained target-tracking scenario. The ablation study focuses on the offline PINN-based solution of the FKE. In this study, all components of the filtering framework remain unchanged except for the computational domain strategy.

\subsubsection{Fixed-domain baseline construction}

In the proposed method, the local computational domain is
\begin{equation}
\Omega_i=[\bar{s}_{i-1}^*-4,\bar{s}_{i-1}^*+4]\times[\bar{v}_{i-1}^*-4,\bar{v}_{i-1}^*+4].
\label{eq:local_domain}
\end{equation}

Monte Carlo simulations show that, over the entire filtering horizon, the trajectory of the density center remains approximately within $s\in[-3,3]$ and $v\in[-3,3]$. To avoid truncation of the density at any timestep, the fixed-domain baseline must cover both the maximal displacement of the density center and the local support width required to resolve the dominant probability mass. Accordingly, the fixed-domain baseline is conservatively defined as
\begin{equation}
\Omega_{\mathrm{fixed}}=[-7,7]\times[-7,7],
\label{eq:fixed_domain}
\end{equation}

Which safely encloses the density support throughout the simulation and covers approximately three times the area of the recentered domain.

\subsubsection{Ablation experiment results}

For a fair comparison, all settings except the domain placement and the resulting grid size are kept identical, including the PINN architecture, loss functions, convergence thresholds, maximum training epochs, and spatial resolution. To maintain the same spatial resolution, the recentered-domain setting uses a $200\times200$ grid with spacing 0.04, whereas the fixed-domain setting uses a $350\times350$ grid to preserve the same spacing over the larger domain. As a result, the number of spatial collocation points increases by approximately a factor of three in the fixed-domain setting.

Fig.~\ref{fig:ablation} illustrates the state trajectories reconstructed from the offline FKE solutions under the two domain strategies. The black curve denotes the ground truth, the red curve corresponds to the proposed recentered-domain strategy, and the gray dashed curve represents the fixed-domain baseline. The fixed-domain solution exhibits noticeable deviations and oscillatory behavior, particularly in regions where the state undergoes rapid variation. In contrast, the recentered-domain solution closely follows the true trajectory across the time horizon and remains stable during highly nonlinear phases. This comparison indicates that domain recentering improves the density approximation in the high-probability region, thereby supporting more accurate state recovery in this ablation setting.

\begin{inlinefigure}
\centering
\includegraphics[width=0.8\textwidth]{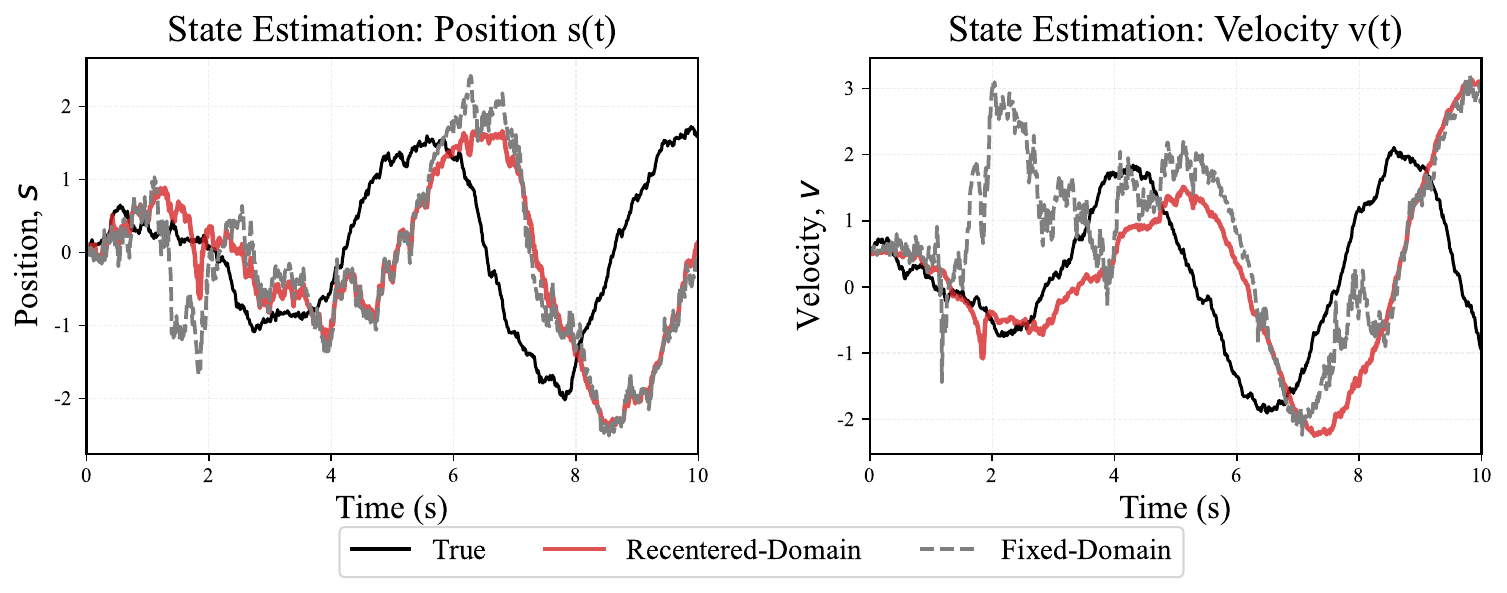}
\caption{Comparison of state estimation results under recentered-domain and fixed-domain strategies.}
\label{fig:ablation}
\end{inlinefigure}

Table~\ref{tab:ablation} summarizes the quantitative comparison between the recentered-domain and fixed-domain strategies. Under the same network architecture and uniform collocation strategy, the fixed-domain baseline suffers from reduced effective resolution in the high-probability region. In addition, the average number of training epochs per timestep decreases from 8436 to 3568 under the recentered-domain setting, indicating faster convergence of the PINN solver.

\begin{inlinetable}
\caption{Comparison between recentered-domain and fixed-domain strategies in terms of estimation accuracy and training convergence.}
\label{tab:ablation}
\centering
\begin{tabular}{@{}lccc@{}}
\toprule
Domain strategy & MSE-$s$ & MSE-$v$ & Avg. epochs per timestep\\
\midrule
Fixed domain & 0.559 & 1.327 & 8436\\
Recentered domain & 0.145 & 0.179 & 3568\\
\bottomrule
\end{tabular}
\end{inlinetable}

The ablation results indicate that the performance degradation observed in the fixed-domain setting is associated with reduced effective numerical resolution within the high-probability region, rather than only with increased computational cost. When a large static domain is used, the dominant probability mass occupies only a small portion of the computational region. Consequently, collocation points are distributed over extensive low-probability areas, reducing the effective approximation capacity of the PINN in the region that critically influences state estimation. As a result, both density reconstruction quality and filtering accuracy deteriorate.

In contrast, the recentered-domain strategy maintains a fixed-size computational window that continuously tracks the evolving high-probability region. This concentrates numerical resolution where it is most needed and preserves approximation fidelity without increasing network capacity. The improvement in estimation accuracy, together with the reduction in average training epochs, indicates that domain recentering is a key component of the proposed framework rather than a minor implementation detail.

\FloatBarrier

\section{Conclusion}\label{sec:conclusion}

In this study, we proposed RD-YYF, an offline-online decoupled framework for efficient state estimation in nonlinear filtering problems. RD-YYF introduces a dynamically recentered local computational window to reduce the computational burden of solving the FKE over a large global domain. The offline stage uses PINNs to generate FKE solution snapshots and applies PCA to construct a low-dimensional representation of the density evolution. A lightweight residual surrogate is then trained to predict terminal solution coefficients from initial-condition coefficients and the domain center. During online execution, the pretrained surrogate is deployed within the recentered window to support fast per-timestep density prediction, observation update, and state estimation.

Numerical experiments in two continuous-time target-tracking scenarios showed that RD-YYF achieved lower average estimation errors than EKF and PF under the tested settings. The online timing results further indicate that RD-YYF maintains practical per-timestep inference cost after offline training, although it is not the fastest method among the compared baselines. The ablation study further confirmed the role of the recentered-domain strategy: by concentrating numerical resolution near the high-probability region, the strategy improved density approximation quality and reduced the number of PINN training epochs required in the offline stage.

Future work will extend RD-YYF to higher-dimensional nonlinear filtering problems and investigate its robustness under time-varying noise statistics, uncertain motion constraints, and more complex sensing configurations. These studies will further assess the adaptability of the proposed framework in broader target-tracking scenarios.

\section*{Funding}
This research did not receive any specific grant from funding agencies in the public, commercial, or not-for-profit sectors.

\printcredits

\section*{Declaration of competing interest}
The authors declare that they have no known competing financial interests or personal relationships that could have appeared to influence the work reported in this paper.
\section*{Declaration of generative AI and AI-assisted technologies in the manuscript preparation process}

During the preparation of this work, the authors used ChatGPT to improve
the readability and language of the manuscript. After using this tool,
the authors reviewed and edited the content as needed and take full
responsibility for the content of the publication.

\section*{Data availability}
The processed numerical data used to generate Tables~\ref{tab:ex1} and~\ref{tab:ex2} are provided as supplementary data in the file Ex1\_Ex2\_MSE.xlsx. This file includes a README sheet, offline and online seed lists, per-seed MSE values, and per-timestep online timing results. The seed-data files for all experiments are available in the supporting data repository at https://gitee.com/malei0506/rdyyf.


\begin{thebibliography}{31}
\bibitem{ref1} Akyon, F.C., Eryuksel, O., Ozfuttu, K.A., Altinuc, S.O.: Track boosting and synthetic data aided drone detection. In: 2021 17th IEEE International Conference on Advanced Video and Signal Based Surveillance (AVSS), pp. 1--5. IEEE (2021)
\bibitem{ref2} Ameen, S., Vokhidov, H.: Autonomous mobile robot navigation: Tracking problem (2024)
\bibitem{ref3} Zheng, D., Sun, Y., He, D., Wang, J., Ju, H., Zhao, H.: Application of UAV target detection and tracking system based on intelligent image recognition. In: 2024 2nd World Conference on Communication \& Computing (WCONF) (2024)
\bibitem{ref4} Kalman, R.E.: A new approach to linear filtering and prediction problems. Journal of Basic Engineering 82(1), 35--45 (1960)
\bibitem{ref5} Chen, H., Feng, X., Huang, Z., Zhuang, Z.: Target tracking based on Kalman filtering techniques. In: 2022 International Symposium on Control Engineering and Robotics (ISCER), pp. 237--245 (2022)
\bibitem{ref6} Auger, F., Hilairet, M., Guerrero, J.M., Monmasson, E., Orlowska-Kowalska, T., Katsura, S.: Industrial applications of the Kalman filter: A review. IEEE Transactions on Industrial Electronics 60(12), 5458--5471 (2013)
\bibitem{ref7} Gordon, N.J., Salmond, D.J., Smith, A.F.M.: Novel approach to nonlinear/non-Gaussian Bayesian state estimation. IEE Proceedings F 140(2), 107--113 (1993)
\bibitem{ref8} Zhang, C., Yang, Y., Ding, Y.: Research on underwater target tracking based on Gaussian Hermitian Kalman particle filter algorithm. Journal of Physics: Conference Series (2021)
\bibitem{ref9} Li, G., Qiao, Y.: A ship target detection and tracking algorithm based on graph matching. Journal of Physics: Conference Series (2021)
\bibitem{ref10} Li, Y., Lou, J., Tan, X., Xu, Y., Zhang, J., Jing, Z.: Adaptive kernel learning Kalman filtering with application to model-free maneuvering target tracking. IEEE Access (2022)
\bibitem{ref11} Aranda, I.A., Pérez-Zúñiga, G.: Highly maneuverable target tracking under glint noise via uniform robust exact filtering differentiator with intrapulse median filter. IEEE Transactions on Aerospace and Electronic Systems (2022)
\bibitem{ref12} Revach, G., Shlezinger, N., Ni, X., Escoriza, A.L., van Sloun, R.J.G., Eldar, Y.C.: KalmanNet: Neural network aided Kalman filtering for partially known dynamics. IEEE Transactions on Signal Processing 70, 1532--1547 (2022)
\bibitem{ref13} Song, F., Li, Y., Cheng, W., Dong, L., Li, M., Li, J.: An improved Kalman filter based on long short-memory recurrent neural network for nonlinear radar target tracking. Wireless Communications and Mobile Computing (2022)
\bibitem{ref14} Shen, L., Su, H., Li, Z., Jia, C., Yang, R.: Self-attention-based transformer for nonlinear maneuvering target tracking. IEEE Transactions on Geoscience and Remote Sensing (2023)
\bibitem{ref15} Das, S., Kumar, K., Bhaumik, S.: Tracking an underwater target with unknown measurement noise statistics using variational Bayesian filters. arXiv: Electrical Engineering and Systems Science, Signal Processing (2023)
\bibitem{ref16} Yau, S.-T., Yau, S.S.-T.: Real-time solution of the nonlinear filtering problem without memory I. Mathematical Research Letters 7(6), 671--693 (2000)
\bibitem{ref17} Yau, S.-T., Yau, S.S.-T.: Real-time solution of the nonlinear filtering problem without memory II. SIAM Journal on Control and Optimization 47(1), 163--195 (2008)
\bibitem{ref18} Duncan, T.E.: Probability densities for diffusion processes with applications to nonlinear filtering theory and detection theory. Stanford University (1967)
\bibitem{ref19} Mortensen, R.E.: Optimal control of continuous-time stochastic systems. University of California, Berkeley (1966)
\bibitem{ref20} Zakai, M.: On the optimal filtering of diffusion processes. Zeitschrift für Wahrscheinlichkeitstheorie und verwandte Gebiete 11, 230--243 (1969)
\bibitem{ref21} Shi, J., Yang, Z., Yau, S.S.-T.: Direct method for Yau filtering system with nonlinear observations. International Journal of Control 91(3), 678--687 (2018)
\bibitem{ref22} Chen, X., Shi, J., Yau, S.S.-T.: Real-time solution of time-varying Yau filtering problems via direct method and Gaussian approximation. IEEE Transactions on Automatic Control 64(4), 1648--1654 (2019)
\bibitem{ref23} Raissi, M., Perdikaris, P., Karniadakis, G.E.: Physics-informed neural networks: A deep learning framework for solving forward and inverse problems involving nonlinear partial differential equations. Journal of Computational Physics 378, 686--707 (2019)
\bibitem{ref24} Davis, M.H.A., Marcus, S.I.: An introduction to nonlinear filtering. In: Hazewinkel, M., Williams, J.S. (eds.) Stochastic Systems: The Mathematics of Filtering and Identification and Applications, pp. 53--75. Reidel, Dordrecht (1981)
\bibitem{ref25} Karniadakis, G.E., Kevrekidis, I.G., Lu, L., Perdikaris, P., Wang, S., Yang, L.: Physics-informed machine learning. Nature Reviews Physics 3, 422--440 (2021)
\bibitem{ref26} Abdi, H., Williams, L.J.: Principal component analysis. Wiley Interdisciplinary Reviews: Computational Statistics 2(4), 433--459 (2010)
\bibitem{ref27} Benner, P., Gugercin, S., Willcox, K.: A survey of projection-based model reduction methods for parametric dynamical systems. SIAM Review 57(4), 483--531 (2015)
\bibitem{ref28} Rao, C.V., Rawlings, J.B., Mayne, D.Q.: Constrained state estimation for nonlinear discrete-time systems: stability and moving horizon approximations. IEEE Transactions on Automatic Control 48(2), 246--258 (2003)
\bibitem{ref29} Ali, T.M.F.: Maneuvering, multi-target tracking using particle filters. arXiv:1910.04379 (2019)
\bibitem{ref30} Kingma, D.P., Ba, J.: Adam: A method for stochastic optimization. In: International Conference on Learning Representations (2015)
\bibitem{ref31} He, K., Zhang, X., Ren, S., Sun, J.: Deep residual learning for image recognition. In: Proceedings of the IEEE Conference on Computer Vision and Pattern Recognition, pp. 770--778 (2016)
\end{thebibliography}
\end{document}